\documentclass[11pt]{amsart}
\usepackage {amsmath}
\usepackage {amsthm}
\usepackage {enumerate}
\usepackage{amssymb}
\usepackage{amsfonts}
\usepackage{anysize}
\usepackage{setspace}
\usepackage[all]{xy}

\marginsize{3.15cm}{3.15cm}{1.6cm}{3cm}

\newcommand{\N}{\mathbb N}

\newcommand{\R}{\mathbb R}

\numberwithin{equation}{section}

\theoremstyle{plain}

\newtheorem{theorem}[equation]{Theorem}
\newtheorem{proposition}[equation]{Proposition}
\newtheorem{corollary}[equation]{Corollary}
\newtheorem{lemma}[equation]{Lemma}
\newtheorem{question}[equation]{Question}

\theoremstyle{definition}
\newtheorem{definition}[equation]{Definition}

\theoremstyle{definition}
\newtheorem{remark}[equation]{Remark}

\theoremstyle{definition}

\theoremstyle{definition}

\theoremstyle{definition}

\title{\scshape\bfseries Characterizing Corson and Valdivia compact spaces}

\author{{\bfseries F. Casarrubias-Segura}}

\address{Departamento de Matematicas\\
         Facultad de Ciencias\\
         UNAM\\
         Circuito exterior s/n, Ciudad Universitaria, C.P. 04510, M\'exico D. F., Mexico}
\email{fcasarrubiass@ciencias.unam.mx}

\author{{\bfseries S. Garc\'ia-Ferreira}}
\address{Centro de Ciencias Matem\'aticas\\
         Universidad Nacional Aut\'onoma de M\'exico\\
				 Campus Morelia\\
         Apartado Postal 61-3, Santa Mar\'ia, 58089, Morelia, Michoac\'an, M\'exico.}
\email{sgarcia@matmor.unam.mx}

\author{{\bfseries R. Rojas-Hern\'andez}}
\email{satzchen@yahoo.com.mx}

\subjclass[2010]{Primary 54C50, 54C15, 54E45.}

\keywords{Corson compact, Valdivia compact, $r$-skeleton, $q$-skeleton, $c$-skeleton, monotonically retractable, function space}

\date{}

\thanks{Research of the first-named author was supported
by  CONACYT grant no. 176202 and PAPIIT grant no. IN-101911.}

\begin{document}

\begin{abstract} We give a new characterization of Valdivia compact spaces: A compact space is Valdivia if and only if it has a dense commutatively monotonically retractable subspace.
This result solves Problem 5.12 from \cite{sal-rey}. Besides,  we introduce the notion of full $c$-skeleton and prove that a compact space is Corson if and only if it has a full $c$-skeleton.
\end{abstract}

\maketitle

\section{Introduction}

Our spaces will be Tychonoff (completely regular and Hausdorff). The Greek letter $\omega$ will stand for the first infinite cardinal number. Given an infinite set  $X$, $\mathcal{P}(X)$ is the power set of $X$, the symbol $[X]^{\leq \omega}$ will denote the set of all countable  subsets of $X$ and the meaning of $[X]^{< \omega}$ should be clear.  The  real line with the usual order topology will be denoted by  $\R$ and $\mathcal{B}(\R)$ will stand for a countable fixed base for the topology of  $\R$. If $T$ is any set we set $\Sigma\R^T = \{f \in \R^T : \left|f^{-1}(\R \setminus \{0\})\right| \leq \omega\}$. For a space $X$, $C_p(X)$ will be the set $C(X)$ of all real-valued continuous functions on $X$ equipped with the topology of pointwise convergence and $\mathcal{CL}(X)$ will be the family of all non-empty closed subsets of $X$.
For a continuous map $f : X \to Y$ we denote by $f^{\ast} : C_p(Y) \to C_p(X)$ the dual map of $f$ given by $f^{\ast}(g) = g \circ f$ for all $g \in C_p(Y)$. If $Y \subseteq X$, then we denote by $\pi_Y : C_p(X) \to C_p(Y)$ the function which restricts any map in $C_p(X)$ to $Y$.
For a nonempty set $A \subseteq C_p(X)$ the map $\Delta_{A} : X \to \R^A$ is called the {\it diagonal map} of $A$.
A onto map $f : X \to Y$ is called {\it ${\mathbb R}$-quotient} if, for every function  $g : Y \to {\mathbb R}$ the continuity of the composition $g \circ f$ implies the continuity of $g$.
A space $X$ is called {\it cosmic} if it has a countable network. For  $\mathcal{N} \subseteq \mathcal{P}(X)$ and $f : X \to Y$, we say that $\mathcal{N}$ is a {\it network of $f$} if for every $x \in X$ and each open set $U$ in $Y$ with $f(x) \in U$ there is $N \in \mathcal{N}$ such that $x \in  N$ and $f (N) \subseteq U$.

\medskip

In this paper $\Gamma$ always will denote an up-directed $\sigma$-complete partially ordered set. Given $\Gamma$ and a set $Y$, a function  $\phi : \Gamma \to [Y]^{\leq\omega}$  is called $\omega$-{\it monotone} provided that:
\begin{enumerate}[(a)]
	\item if $s,t \in \Gamma$ and $s \leq t$ then $\phi(s) \subseteq \phi(t)$; and
  \item if $\{s_n : n \in \mathbb{N}\} \subseteq \Gamma$ and $s_n \leq s_{n+1}$ for each $n \in \mathbb{N}$, then $\phi(\sup_{n \in \mathbb{N}} s_n) =  \bigcup\{\phi(s_n): n \in \N\}.$
\end{enumerate}
For the case when $\Gamma = [X]^{\leq\omega}$, where $X$ is a nonempty set, the order involved will be always  the containment. Given an $\omega$-monotone map $\phi : [X]^{\leq\omega} \to [X]^{\leq\omega}$, we consider the $\omega$-monotone map $\overline{\phi} : [X]^{\leq\omega} \to [X]^{\leq\omega}$ defined by
\begin{center}
$\overline{\phi}(A) = \bigcap\left\{B \in [X]^{\leq\omega} : A \subseteq B \ \text{and} \ \phi(B) \subseteq B\right\},$
\end{center}
for all $A \in [X]^{\leq\omega}$.
All topological notions whose definitions are not stated explicitly here should be understood as  in  \cite{ark} and \cite{ltka} .

\begin{definition}
A set $Y \subseteq X$ will be called a \textit{$\Sigma$-subset} of $X$ if there is an embedding $\phi : X \to \R^T$, for some set $T$, such that $Y = \phi^{-1}(\Sigma \R^T)$. A compact space is called \textit{Valdivia} if it admits a dense $\Sigma$-subset. A compact space is \textit{Corson} if it is a $\Sigma$-subset of itself.
\end{definition}

It is well known that  Valdivia compact spaces play an important role in the study of the structure of nonseparable Banach spaces.
 Having certain nonseparable Banach space, sometimes it is useful to decompose it into smaller pieces (subspaces).  One possible  decomposition is the notion of a projectional resolution of the identity (see e.g. \cite[Definition 3.35]{HMVZ}).  Another and very useful decomposition is the concept of projectional skeleton introduced by W.Kubis in \cite{kbs}. The dual notion of projectional skeleton is the following.

\begin{definition}{\bf\cite{kbs-mich}}\label{DRS}
An {\it $r$-skeleton} in a space $X$ is a family $\{r_s : s \in \Gamma\}$ of retractions satisfying:
\begin{enumerate}[$(i)$]
	\item  $r_s(X)$ is cosmic for each $s \in \Gamma$;
  \item  $r_s = r_s \circ r_t = r_t \circ r_s$ whenever $s \leq t$;
  \item if $\{s_n : n \in \mathbb{N}\} \subseteq \Gamma$, $s_n \leq s_{n+1}$ for each $n \in \mathbb{N}$ and $t = \sup_{n \in \mathbb{N}} s_n$, then $r_t(x) =  \lim_{n \to \infty} r_{s_n}(x)$ for each $x \in X$; and
  \item $x = \lim_{s \in \Gamma} r_s(x)$ for every $x \in X$.
\end{enumerate}
We shall say that $\bigcup\{r_s(X) : s \in \Gamma\}$ is the subspace induced by the $r$-skeleton $\{r_s : s \in \Gamma\}$. When $X = \bigcup\{r_s(X) : s \in \Gamma\}$, we say that  the $r$-skeleton $\{r_s : s \in \Gamma\}$ is  {\it full}. The $r$-skeleton $\{r_s : s \in \Gamma\}$ is called \textit{commutative} if $r_s \circ r_t = r_t \circ r_s$ for every $s, t \in \Gamma$.
\end{definition}

A characterization of the class of all Valdivia compacta in terms of a system of retractions was obtained in \cite{kbs-mich}. Certainly, the authors proved that a compact space is Valdivia if and only if it has a commutative $r$-skeleton.  Our of main purpose in the first section is to apply a method similar to the one used in \cite{ban2} to give another proof of this  characterization of Valdivia compact spaces. We use this characterization to answer a question in \cite{sal-rey}, by   characterizing   Valdivia compact spaces in terms of families of networks and retractions.
We know that Corson compact spaces are characterized as compact spaces with a full $r$-skeleton \cite{mrk2}.
In \cite{sal-rey}, function spaces over Corson compact spaces were characterized by using families of $\R$-quotient maps and countable sets (full $q$-skeletons), in an analogous sense to $r$-skeletons. In the last part of this paper, we provide a new characterization of Corson compact spaces by using  families of closed sets and continuous maps ($c$-skeletons).

\section{A characterization of Valdivia compact spaces}

In order to give new proof that every compact space with a commutative $r$-skeleton is Valdivia, we establish some  basic facts.

\begin{lemma}\label{LL}
Let $\{r_s : s \in \Gamma\}$ be a family of retractions on a countably compact space $X$ satisfying $(i)$ - $(iii)$ of Definition \ref{DRS}. If $Y = \bigcup\{r_s(X) : s \in \Gamma\}$, then for each $\mathcal{\mathcal{A}} \in [\mathcal{P}(Y)]^{\leq\omega}$ and $s_0 \in \Gamma$ there exist $s \in \Gamma$ and $D \in [\bigcup\mathcal{A}]^{\leq\omega}$ such that $s_0 \leq s$ and $r_s(\textnormal{cl}(A)) = \textnormal{cl}(D \cap A)$ for all $A \in \mathcal{A}$.
\end{lemma}

\proof
For each $n \in \N$, we will construct $D_n \in [\bigcup\mathcal{A}]^{\leq\omega}$ and $s_{n+1} \in \Gamma$ recursively as follows. Choose $D_{n} \in [\bigcup\mathcal{A}]^{\leq\omega}$ such that, for each $A \in \mathcal{A}$, the set $r_{s_n}(D_n \cap A)$ is dense in $r_{s_n}(\textnormal{cl}(A))$. Besides, select $s_{n+1} \in \Gamma$ such that $s_{n} \leq s_{n + 1}$ and $D_n \subseteq r_{s_{n + 1}}(X)$. Let $D = \bigcup\{D_n : n \in \N\}$ and $s = \sup\{s_n : n \in \N\}$. We shall verify that $D$ and $s$ are the required sets. Choose $A \in \mathcal{A}$ and $x \in \textnormal{cl}(A)$. Let $U$ be a neighborhood of $r_s(x)$. Fix an open set $V$ in $X$ such that $r_s(x) \in V \subseteq \textnormal{cl}(V) \subseteq U$. Then we can find $N \in \N$ such that $r_{s_n}(x) \in V$ whenever $n \geq N$. For each $n \geq N$, because of $(ii)$, there is $y_n \in D_n \cap A$ such that $r_{s_m}(y_n) \in V$ whenever $N \leq m \leq n$. For every $n < N$ select $y_n \in D_n \cap A$ arbitrarily. Observe that $\{y_n : n \in \N\} \subseteq r_s(X)$. Since $r_s(X)$ is compact, there is an accumulation point $y$ of $\{y_n : n \in \N\}$ in $r_s(X)$.  We can assume that $y = \lim_{n \to \infty} y_n$. Then $r_s(y) = \lim_{m \to \infty} r_{s_m}(y) = \lim_{m \to \infty} \lim_{n \to \infty} r_{s_m}(y_n) \in \textnormal{cl}(V) \subseteq U$. On the other hand, we know that  $r_s(y) = \lim_{n \to \infty} r_s(y_n) = \lim_{n \to \infty} y_n$, and hence we must have that $y_n \in U \cap (D \cap A)$ for $n$ large enough. Therefore, $r_s(x) \in \textnormal{cl}(D \cap A)$.
\endproof

\begin{corollary}\label{LTL}
Let $\{r_s : s \in \Gamma\}$ be a family of retractions in a countably compact space $X$ satisfying $(i)$ - $(iii)$ of Definition \ref{DRS}. If $Y = \bigcup\{r_s(X) : s \in \Gamma\}$, then
\begin{enumerate}[(1)]
	\item $t(Y) \leq \omega$.
	\item $x = \lim_{s \in \Gamma}r_s(x)$ for each $x \in \textnormal{cl}(Y)$.
\end{enumerate}
\end{corollary}

\proof $(1)$ Assume that $A \subseteq Y$ and $x \in \textnormal{cl}_Y(A)$. Set $\mathcal{A} = \{A\}$ and choose $s_0 \in \Gamma$ so that $x \in r_{s_0}(X)$. By Lemma \ref{LL}, we can find $D \in [A]^{\leq\omega}$ and $s \in \Gamma$ for which $s_0 \leq s$ and $r_s(\textnormal{cl}(A)) = \textnormal{cl}(D)$. This implies that $x = r_s(x) \in r_s(\textnormal{cl}(A)) = \textnormal{cl}(D)$.

\medskip

$(2)$ Fix $x \in \textnormal{cl}(Y)$. Let $U$ be an open neighborhood of $x$. Choose an open set $V$ in $X$ such that $x \in V \subseteq \textnormal{cl}(V)\subseteq U$. Set $A_1 = V \cap Y$, $A_2 = (X \setminus U) \cap Y$ and $\mathcal{A} = \{A_1,A_2\}$. By Lemma \ref{LL}, we can find $s \in \Gamma$ such that $r_s(\textnormal{cl}(A_i)) \subseteq \textnormal{cl}(A_i)$ for $i = 1,2$. As $x \in \textnormal{cl}(A_1)$, then  $r_s(x) \in \textnormal{cl}(A_1)$. Choose $t \in \Gamma$ such that $s \leq t$. Assume that  $r_t(x) \not \in U$. Then $r_t(x) \in (X \setminus U) \cap Y = A_2$ and so $r_s(x) = r_s(r_t(x)) \in \textnormal{cl}(A_2)$, but this a contradiction since $\textnormal{cl}(A_1) \cap \textnormal{cl}(A_2) = \emptyset$. Thus, $r_t(x) \in U$.
\endproof

\begin{lemma}\label{LDH}
Let $X$ be compact and let $C$ be closed in $X$. Assume that $\{r_s : s \in \Gamma\}$ is a family of retractions from $X$ into $C$ such that $\{r_s\restriction C : s \in \Gamma\}$ is an $r$-skeleton on $C$. If $R = \Delta_{\{r_s : s \in \Gamma\}}$, then $R \restriction_C : C \to R(X)$ is a homeomorphism.
\end{lemma}

\proof
The map $R\restriction_C$ is clearly continuous. From condition $(iv)$ we deduce that $R \restriction C$ is one-to-one. Since $C$ is compact, we only need to verify that $R(C) = R(X)$. Set $Y = \bigcup\{r_s(X) : s \in \Gamma\} \subseteq C$. We claim that $R(Y)$ is dense in $R(X)$. Indeed, choose a nonempty open set $U$ of $R(X)$. We can assume that there exists a finite set $F \subseteq \Gamma$ and $U_s \subseteq r_s(X)$, for every $s \in F$,  so that $U = \{R(x) : x \in X \textnormal{ and } r_s(x) \in U_s \textnormal{ for each }s \in F\}$. Select $x \in X$ such that $R(x) \in U$. Let $t \in \Gamma$ be an upper bound of $F$. Note that $r_t(x) \in Y$. Moreover, since $R(x) \in U$, we deduce that $r_s(r_t(x)) = r_s(x) \in U_s$ for each $s \in F$; that is, $R(r_t(x)) \in U$. Therefore, $R(Y)$ is dense in $R(X)$.  The compactness of $C$ implies that $R(C) = R(X)$.
\endproof

 The next lemma was used in \cite{sal-rey} to construct, from a given $r$-skeleton, another $r$-skeleton indexed by the  family of countable subsets of the space on which it is defined. For our purposes, it will be useful to get a suitable structure of retractions indexed by the family of all subsets of a given subspace.

\begin{lemma}\label{LAMG}
Given  a set $X$ and $\Gamma$, suppose that for each $x \in X$ we have assigned $s_x \in \Gamma$. Then there exists a function  $\gamma : [X]^{\leq\omega} \to \Gamma$ such that:
\begin{enumerate}[(1)]
  \item $\gamma(\{x\}) \geq s_x$ for each $x \in X$.
	\item If $A, B \in [X]^{\leq\omega}$ and $A \subseteq B$, then $\gamma(A) \leq \gamma(B)$.
	\item If $\{A_n : n \in \mathbb{N}\} \subseteq [X]^{\leq\omega}$ and $A_n \leq A_{n+1}$ for each $n \in \mathbb{N}$, then

\centerline{$\gamma(\bigcup\{A_n : n \in \mathbb{N}\}) =  \sup_{n \in \mathbb{N}} \gamma(A_n)$.}
\end{enumerate}
\end{lemma}

\begin{lemma}\label{LRA}
Let $X$ be compact and let $Y$ be induced by a commutative $r$-skeleton $\{r_s : s \in \Gamma\}$ in $X$. Then there exists a family $\{r_A : A \in \mathcal{P}(Y)\}$ of retractions in $X$ such that for every $A \in \mathcal{P}(Y)$ the following conditions hold:
\begin{enumerate}[(a)]
    \item $A \subseteq r_A(X)$ and $d(r_A(X)) \leq \left|A\right|$.
	\item The family $\{r_B \restriction_{r_A(X)} : B \in [A]^{\leq\omega}\}$ is a commutative $r$-skeleton on $r_A(X)$ and induces $Y \cap r_A(X)$.
	\item $r_B \circ r_A = r_A \circ r_B = r_B$ whenever $B \subseteq A$.
    \item $r_A(Y) \subseteq Y$.
\end{enumerate}
\end{lemma}

\proof First, we define $r_A $ for each $A \in [Y]^{\leq \omega}$. Given $y \in Y$ fix $s_y \in \Gamma$ such that $y \in r_{s_y}(X)$. Then, we consider $\gamma : [Y]^{\leq\omega} \to \Gamma$ as in the Lemma \ref{LAMG}. Thus, for every $A \in [Y]^{\leq\omega}$ we define $r_A = r_{\gamma(A)}$. Now, we shall proceed to define $r_A$ for an uncountable $A \subseteq Y$:
For each $F \in [Y]^{<\omega}$ fix a countable dense subset $D_F$ of $r_F(X)$ such that $F \subseteq D_F$. Consider the map $\mathcal{D} : \mathcal{P}(Y) \to \mathcal{P}(Y)$ given by $\mathcal{D}(A) = \bigcup\{D_F : F \in [A]^{<\omega}\}$ for all $A \in \mathcal{P}(Y)$. Note that $A \subseteq \mathcal{D}(A)$ and $\left|\mathcal{D}(A)\right| \leq \left|A\right|$ for each $A \in \mathcal{P}(Y)$. We remark  that  $r_A(X) = \textnormal{cl}(\mathcal{D}(A))$ whenever $A \in [Y]^{\leq\omega}$. Fix $A \in \mathcal{P}(Y) \setminus [Y]^{\leq\omega}$. We assert that $Y \cap \textnormal{cl}(\mathcal{D}(A)) = \bigcup\{r_B(X) : B \in [A]^{\leq\omega}\}$. In one hand, take $y \in Y \cap \textnormal{cl}(\mathcal{D}(A))$. From Corollary \ref{LTL} $(1)$  we can find $D \in [\mathcal{D}(A)]^{\leq\omega}$ such that $y \in \textnormal{cl}(D)$. Choose $B \in [A]^{\leq\omega}$ such that $D \subseteq \mathcal{D}(B)$. Then $y \in \textnormal{cl}(D) \subseteq \textnormal{cl}(\mathcal{D}(B)) = r_B(X)$. On the other hand, choose $x \in r_B(X)$ for some $B \in [A]^{\leq\omega}$. Then $x \in r_B(X) = \textnormal{cl}(\mathcal{D}(B)) \subseteq Y \cap \textnormal{cl}(\mathcal{D}(A))$.
This shows the assertion. Thus, we have that  $\{r_B  : B \in [A]^{\leq\omega}\}$ is a family of retractions from $X$ into $\textnormal{cl}(\mathcal{D}(A))$. It follows from the properties of $\gamma$ and Corollary  \ref{LTL} $(2)$ that $\{r_B \restriction_{\textnormal{cl}(\mathcal{D}(A))} : B \in [A]^{\leq\omega}\}$ is a commutative $r$-skeleton on $\textnormal{cl}(\mathcal{D}(A))$. Then, by applying Lemma \ref{LDH}, we can see that  $R_A \restriction \textnormal{cl}(\mathcal{D}(A)) : \textnormal{cl}(\mathcal{D}(A)) \to R_A(X)$ is a homeomorphism where $R_A = \Delta_{\{ r_B: B \in [A]^{\leq\omega}\}}$. Finally, define $r_A = (R_A \restriction \textnormal{cl}(\mathcal{D}(A)))^{-1} \circ R_A$.

\medskip

The clauses $(a)$ and $(b)$ follows directly from the construction. The clauses $(c)$ and $(d)$ are direct consequences from the definition  when $A \in [Y]^{\leq\omega}$.  So we only need to  prove these two conditions for  an arbitrary $A \in \mathcal{P}(Y) \setminus [Y]^{\leq\omega}$.

\medskip

$(c)$ Fix $B \subseteq A$. By construction, we know that  $r_B(X) = \textnormal{cl}(\mathcal{D}(B))  \subseteq  \textnormal{cl}(\mathcal{D}(A))  = r_A(X)$, and so $r_A \circ r_B = r_B$. To verify that $r_B = r_B \circ r_A$ we pick $x \in X$. Note that $R_A(x) = R_A(r_A(x))$; i.e., $r_{D}(x) = r_{D}(r_A(x))$ for each $D \in [A]^{\leq\omega}$. Hence, $r_D(x) = r_D(r_A(x))$ for each $D \in [B]^{\leq\omega}$, and so $R_B(x) = R_B(r_A(x))$. Therefore, $r_B(x) = r_B(r_A(x))$.\smallskip

$(d)$ Fix $y \in Y$ and let $B = \{y\}$. It is enough to show that $r_A(r_B(y)) = r_B(r_A(y))$. Apply $(b)$ to see that $r_A(y) = \lim_{D \in [A]^{\leq\omega}} r_D(r_A(y))$ and, as a consequence,
$$
r_B(r_A(y)) = \lim_{D \in [A]^{\leq\omega}} r_B(r_D(r_A(y))) = \lim_{D \in [A]^{\leq\omega}} r_D(r_B(r_A(y))) \in r_A(X).
$$
By clauses $(b)$ and  $(c)$, we obtain that $r_D(r_B(y)) = r_B(r_D(y)) = r_B(r_D(r_A(y))) = r_D(r_B(r_A(y)))$ for each $D \in [A]^{\leq\omega}$. It then follows that $R_A(r_B(y)) = R_A(r_B(r_A(y)))$. Since $r_B(r_A(y)) \in r_A(X)$ we must have $r_A(r_B(y)) = r_B(r_A(y))$.
\endproof

\begin{theorem}\label{ICRSSC}
Let $X$ be a compact and let $Y$ be a dense subspace of $X$. If $Y$ is induced by a commutative $r$-skeleton on $X$, then $Y$ is a $\Sigma$-subset of $X$.
\end{theorem}

\proof
We will prove the result by induction on the density of $Y$. For $d(Y) = \omega$ the result can be established directly. Assume that $d(Y) = \kappa > \omega$ and the result holds for spaces of density smaller than $\kappa$. Choose the family $\{r_A : A \in \mathcal{P}(Y)\}$ of retractions on $X$ satisfying $(a)$-$(d)$ of Lemma \ref{LRA}. Let $\{y_\alpha : \alpha < \kappa\}$ be a dense subspace of $Y$. For each $\alpha \leq \kappa$, set $A_\alpha = \{y_\beta : \beta < \alpha\}$ and $r_\alpha = r_{A_\alpha}$. Given $\alpha < \kappa$, because of $(a)$ and $(b)$ we may apply the inductive hypothesis to find a set $T_\alpha$ and an embedding $\phi_\alpha : r_\alpha(X) \to \R^{T_\alpha}$ such that $Y \cap r_\alpha(X) = \phi_\alpha^{-1}(\Sigma\R^{T_\alpha})$. We can assume that $T_\alpha \cap T_\beta = \emptyset$ whenever $\alpha \not= \beta$. Consider the set $T = \bigcup_{\alpha < \kappa}T_\alpha$. We identify $\R^T$ with $\prod_{\alpha < \kappa}\R^{T_\alpha}$. Define $\phi : X \to \R^T$ as follows: If $x \in X$ and $\alpha < \kappa$, then
$$
\phi(x)(\alpha) =
\begin{cases}
\phi_{\alpha+1}(r_{\alpha + 1}(x)) - \phi_{\alpha + 1}(r_{\alpha}(x)) &\mbox{if } \alpha > 0; \\
\phi_{0}(r_{0}(x))                                                    &\mbox{if } \alpha = 0.
\end{cases}
$$

To see that $\phi$ is an embedding we only need to show that $\phi$ is one-to-one. Fix  distinct points $x, y \in X$. It follows from $(b)$ applied to $A_\kappa$ together with $(iv)$ and $(iii)$ that there exist $F \in [A_\kappa]^{<\omega}$ such that $r_F(x) \not= r_F(y)$. As a consequence $\{\alpha < \kappa : r_{\alpha}(x) \not= r_{\alpha}(y)\} \not= \emptyset$. Set $\beta = \min\{\alpha < \kappa : r_{\alpha}(x) \not= r_{\alpha}(y)\}$. If $\beta = 0$ it is clear that $\phi(x) \not= \phi(y)$. Otherwise, note that $(b)$ applied to $A_\beta$ together with $(iv)$ and $(iii)$ imply that $\beta = \alpha + 1$ for some ordinal $\alpha < \kappa$. So we deduce that $\phi(x)(\alpha) \not = \phi(y)(\alpha)$, and hence $\phi(x) \not= \phi(y)$.

To prove that $\phi(Y) \subseteq \Sigma \R^T$ it is enough to show that $\{r_\alpha(x) : \alpha < \kappa\}$ is countable for every $x \in Y$. Indeed, fix $x \in Y$and assume that   $f : \omega_1 \to \kappa$ is an increasing map. We assert that  there is $\zeta < \omega_1$ such that  $r_{f(\beta)}(x) = r_{f(\zeta)}(x)$ for every $\zeta < \beta < \omega_1$. Let $\xi = \sup\{f(\alpha) : \alpha < \omega_1\}$. Let us see that $r_\xi(x) \in \textnormal{cl}(\{r_{f(\alpha)}(x) : \alpha < \omega_1\})$. Let $U$ be a neighborhood of $r_\xi(x)$. Pick an open set $V$ such that $r_\xi(x) \in V \subseteq \textnormal{cl}(V) \subseteq U$. The equality $r_\xi(x) = \lim_{B \in [A_\xi]^{\leq\omega}} r_B(r_\xi(x)) = \lim_{B \in [A_\xi]^{\leq\omega}} r_B(x)$ implies that there exists $D \in [A_\xi]^{\leq\omega}$ such that $r_B(x) \in V$ whenever $B \in [A_\xi]^{\leq\omega}$ and $D \subseteq B$. Choose $\alpha < \omega_1$ such that $D \subseteq A_{f(\alpha)}$. The equality $r_{f(\alpha)}(x) = \lim_{B \in [A_{f(\alpha)}]^{\leq\omega}} r_B(r_{f(\alpha)}(x)) = \lim_{B \in [A_{f(\alpha)}]^{\leq\omega}} r_B(x)$ implies that $r_{f(\alpha)}(x) \in \textnormal{cl}(\{r_B(x) : B \in [A_{f(\alpha)}]^{\leq\omega} \textnormal{ and } D \subseteq B\}) \subseteq \textnormal{cl}(V) \subseteq U$.  Since $t(Y) \leq \omega$, we can find $\zeta < \omega_1$ such that $r_\xi(x) \in \textnormal{cl}(\{r_{f(\alpha)}(x) : \alpha < \zeta\}) \subseteq r_{f(\zeta)}(X)$. Then, for each $\beta \geq \zeta$, we deduce that $r_{f(\beta)}(x) = r_{f(\beta)}(r_\xi(x)) = r_\xi(x)$.
In this way, we have shown that $\{r_\alpha(x) : \alpha < \kappa\}$ is countable.

It is easy to see that $\phi(X \setminus Y) \cap \Sigma \R^T = \emptyset$. Therefore, $Y$ is a $\Sigma$-subset of $X$.
\endproof

The converse of Theorem \ref{ICRSSC} is known and not hard to prove by using a saturation argument. So, we obtain the next result.

\begin{theorem}\cite{kbs-mich}
A compact space $X$ is Valdivia if and only if  admits a commutative $r$-skeleton.
\end{theorem}

\section{Relative monotonically retractable spaces}

To obtain a characterization of Valdivia compact spaces, which is done in the next section, we introduce a relative version of monotonically retractable spaces which were studied in the paper \cite{rjs1}. We shall prove that, in general, the topological behaviors of a relatively monotonically retractable subset and a $\Sigma$-subset are similar.

\begin{definition}\label{DRMR}
Given a space $X$, a subspace $Y$ of $X$ is \textit{monotonically retractable in} $X$ if we can assign to each $A \in [Y]^{\leq\omega}$ a retraction $r_A : X \to Y$ and a family $\mathcal{O}(A) \in [\mathcal{P}(Y)]^{\leq\omega}$ such that:
\begin{enumerate}[(i)]
\item $A \subseteq r_A(X)$;
\item $\mathcal{O}(A)$ is a network of $r_A \restriction Y$; and
\item $\mathcal{O}$ is $\omega$-monotone.
\end{enumerate}
If, in addition,  $r_A \circ r_B = r_B \circ r_A$ for each $A,B \in [Y]^{\leq\omega}$, then we say that $Y$ is \textit{commutatively monotonically retractable in} $X$.
\end{definition}

Assume that  $Y$ is  monotonically retractable in $X$. It is easy to see that $Y$ is monotonically retractable (in itself). As a consequence, $Y$ is collectionwise normal, $\omega$-stable, $\omega$-monolithic and has countable extent (see Corollary 4.9 from \cite{rjs-tka}). Besides, Corollary $3.20$ of \cite{rjs1} implies that $Y$ has countable tightness. Since $\omega$-monolithic spaces of countable tightness are Fr\'echet-Urysohn, the space $Y$ is Fr\'echet-Urysohn. Finally, it is easy to see that $Y$ is $\omega$-closed.

\medskip

The next lemma will be useful to prove the main properties of relatively monotonically retractable spaces.

\begin{lemma}\label{LDMR}
Assume that $Y$ is monotonically retractable in $X$. Then we can assign to each $\mathcal{F} \in [\mathcal{CL}(Y)]^{\leq\omega}$ a retraction $r_\mathcal{F} : X \to Y$, a set $\mathcal{D}(\mathcal{F}) \in [Y]^{\leq\omega}$ and a family $\mathcal{N}(\mathcal{F}) \in [\mathcal{P}(Y)]^{\leq\omega}$  such that:
\begin{enumerate}[(a)]
\item $\mathcal{D}$ and $\mathcal{N}$ are $\omega$-monotone.
\item $\mathcal{N}(\mathcal{F})$ is a network of $r_\mathcal{F} \restriction Y$ and is closed under finite intersections.
\item If $N \in \mathcal{N}(\mathcal{F}) \cup \{Y\}$, $F \in \mathcal{F} \cup\{Y\}$ and $N \cap F \not= \emptyset$, then $\mathcal{D}(\mathcal{F}) \cap (N \cap F) \not= \emptyset$.
\item $r_\mathcal{F}(X) = \textnormal{cl}(\mathcal{D}(\mathcal{F}))$ and $r_\mathcal{F}(F) \subseteq F$ for each $F \in \mathcal{F}$.
\end{enumerate}
Besides, if  $Y$ is commutatively monotonically retractable in $X$, then  $r_\mathcal{F} \circ r_\mathcal{G} = r_\mathcal{G} \circ r_\mathcal{F}$ for each $\mathcal{F}, \mathcal{G} \in [\mathcal{CL}(Y)]^{\leq\omega}$.
\end{lemma}

\proof
Suppose  that the assignments $A \to r_A$ and $A \to \mathcal{O}(A)$, for each $A \in [Y]^{\leq\omega}$, witness that $Y$ is monotonically retractable in $X$. We can assume that $\mathcal{O}(A)$ is closed under finite intersections for all $A \in [Y]^{\leq\omega}$. For every $N \in \bigcup\{\mathcal{O}(A): A \in [Y]^{\leq\omega}\} \cup \{Y\}$ and $F \in \mathcal{CL}(Y)$ fix a point $y_{N,F} \in N \cap F$ if $N \cap F \not= \emptyset$ and choose $y_{N,F} \in Y$ arbitrarily otherwise.  Define $\mathcal{S} : [\mathcal{CL}(Y)]^{\leq\omega} \to [\mathcal{CL}(Y)]^{\leq\omega}$ by $\mathcal{S}(\mathcal{F}) = \{\{y_{N,F}\} : N \in \mathcal{O}(\{x : \{x\} \in \mathcal{F}\}) \cup \{Y\} \textnormal{ and } F \in \mathcal{F} \cup \{Y\}\}$ for each $\mathcal{F} \in [\mathcal{CL}(Y)]^{\leq \omega}$. It is not hard to prove  that $\mathcal{S}$ is $\omega$-monotone. Consider the $\omega$-monotone map $\overline{S}$. For each $\mathcal{F} \in [\mathcal{CL}(Y)]^{\leq\omega}$ we define $\mathcal{D}(\mathcal{F}) = \{x : \{x\} \in \overline{\mathcal{S}}(\mathcal{F})\}$, $r_\mathcal{F} = r_{\mathcal{D}(\mathcal{F})}$ and $\mathcal{N}(\mathcal{F}) = \mathcal{O}(\mathcal{D}(\mathcal{F}))$. Clauses $(a)$ and $(b)$ are easy to verify. \smallskip

$(c)$ Assume that $N \in \mathcal{N}(\mathcal{F}) \cup \{Y\}$, $F \in \mathcal{F} \cup\{Y\}$ and $N \cap F \not= \emptyset$. Then $N \in \mathcal{O}(\mathcal{D}(\mathcal{F})) \cup \{Y\}$ and, in this way, $\{y_{N,F}\} \in \mathcal{S}(\overline{\mathcal{S}}(\mathcal{F})) = \overline{\mathcal{S}}(\mathcal{F})$. Therefore, $y_{N,F} \in \mathcal{D}(\mathcal{F}) \cap (N \cap F)$.\smallskip

$(d)$ Because of $(c)$ the set $\mathcal{D}(\mathcal{F})$ intersects any member of $\mathcal{N}(\mathcal{F})$. So we can apply $(b)$ to see that $\mathcal{D}(\mathcal{F}) = r_\mathcal{F}(\mathcal{D}(\mathcal{F}))$ is a dense subset of $r_\mathcal{F}(Y) = r_\mathcal{F}(X)$. Now assume that there are $F \in \mathcal{\mathcal{F}}$ and $y \in F$ such that $r_\mathcal{F}(y) \in Y \setminus F$. By $(b)$ we can find $N \in \mathcal{N}(\mathcal{F})$ such that $y \in N$ and $r_\mathcal{F}(N) \subseteq Y \setminus F$. By $(c)$ we can choose $y_{N,F} \in \mathcal{D}(\mathcal{F}) \cap (N \cap F)$. Then $y_{N,F} = r_\mathcal{F}(y_{N,F}) \in r_\mathcal{F}(N) \subseteq Y \setminus F$, but this is impossible. Therefore, $r_\mathcal{F}(F) \subseteq F$ for each $F \in \mathcal{F}$.
\endproof

\begin{proposition}
Let $Y$ be a subspace of $X$. If $Y$ is monotonically retractable in $X$, then $Y$ is $C$-embedded in $X$.
\end{proposition}

\proof
Assign to each $\mathcal{F} \in [\mathcal{CL}(Y)]^{\leq\omega}$ a retraction $r_\mathcal{F} : X \to Y$ as in Lemma \ref{LDMR}. Fix a map $f \in C(Y)$. Choose $\mathcal{F} = \{f^{-1}(\textnormal{cl}(B)) : B \in \mathcal{B}(\R)\} \in [\mathcal{CL}(Y)]^{\leq\omega}$. From clause $(d)$ of Lemma \ref{LDMR}, we know that  $r_\mathcal{F} \restriction_Y (F) \subseteq F$ for each $F \in \mathcal{F}$. It follows from Lemma 4.15 in \cite{rjs-tka} that $f = f \circ r_\mathcal{F} \restriction_Y$. Therefore, $f \circ r_\mathcal{F}$ is the desired extension of $f$.
\endproof

\begin{corollary}
Let $Y$ be a dense subspace of $X$. If $Y$ is monotonically retractable in $X$, then $\beta X = \beta Y$.
\end{corollary}

It is well-known that a dense and $C$-embedded subspace of a given space is $G_\delta$-dense in it (for a proof see \cite[Th. 6.1.4]{ark-tko}). So, we have the following corollary.

\begin{corollary}\label{Gdelta}
Let $Y$ be a dense subspace of $X$. If $Y$ is monotonically retractable in $X$, then $Y$ is $G_\delta$-dense in $X$.
\end{corollary}

\begin{proposition}
Assume that $Y$ is monotonically retractable in $X$. If $\{F_n : n \in \N\} \subseteq \mathcal{CL}(Y)$ and $Z = \textnormal{cl}_X(\bigcup\{F_n : n \in \N\})$, then $Y \cap Z$ is monotonically retractable in $Z$.
\end{proposition}

\proof
Assign to each $\mathcal{F} \in [\mathcal{CL}(Y)]^{\leq\omega}$ a retraction $r_\mathcal{F} : X \to Y$ and a family $\mathcal{N}(\mathcal{F}) \in [\mathcal{P}(Y)]^{\leq\omega}$ as in Lemma \ref{LDMR}. Define $\mathcal{G} : [Y \cap Z]^{\leq\omega} \to [\mathcal{CL}(Y)]^{\leq\omega}$ as $\mathcal{G}(A) = \{\{x\} : x \in A\} \cup \{F_n : n \in \N\}$ for all $A \in [Y \cap Z]^{\leq\omega}$. It is clear that $\mathcal{G}$ is $\omega$-monotone. For each $A \in [Y \cap Z]^{\leq\omega}$ let $s_A = r_{\mathcal{G}(A)}\restriction_Z$ and $\mathcal{O}(A) = \{N \cap Z : N \in \mathcal{N}(\mathcal{G}(A))\}$. We will verify that the assignments $A \to \mathcal{O}(A)$ and $A \to r_{\mathcal{G}(A)}$ witness that $Y \cap Z$ is monotonically retractable in $Z$. First note that $\mathcal{O}$ is $\omega$-monotone. Fix $A \in [Y \cap Z]^{\leq\omega}$. It is clear that $\mathcal{O}(A) \in [\mathcal{P}(Y \cap Z)]^{\leq\omega}$. Besides, the fact that $\mathcal{N}(\mathcal{G}(A))$ is a network of $r_{\mathcal{G}(A)} \restriction_Y$ implies that $\mathcal{O}(A)$ is a network of $s_A \restriction_{(Y \cap Z)}$. By property $(d)$ of Lemma \ref{LDMR}, we have that $s_A(F) \subseteq F$ for each $F \in \mathcal{G}(A)$. In one hand, $s_A(\{x\}) \subseteq \{x\}$ for each $x \in A$, and so $A \subseteq s_A(X)$. On the other hand, $s_A(F_n) \subseteq F_n$ for each $n \in \N$, that is $s_A(\bigcup\{F_n : n \in \N\}) \subseteq \bigcup\{F_n : n \in \N\}$. From the continuity of $s_A$ we deduce that $s_A(Z) \subseteq Z$. Therefore, $s_A$ is a continuous retraction from $Z$ to $Y \cap Z$.
\endproof

\begin{proposition}
Let $Y$ be a dense subspace of $X$. If $Y$ is monotonically retractable in $X$ and $Z = \bigcap\{\textnormal{cl}_X(U_n) : n \in \N\} \neq \emptyset$ where each $U_n$ is open in $X$, then $Y \cap Z$ is dense and monotonically retractable in $Z$.
\end{proposition}

\proof According to Corollary \ref{Gdelta}, we have that $Y \cap Z$ is dense in $Z$.
 Assign to each $\mathcal{F} \in [\mathcal{CL}(Y)]^{\leq\omega}$ a retraction $r_\mathcal{F} : X \to Y$ and a family $\mathcal{N}(\mathcal{F}) \in [\mathcal{P}(Y)]^{\leq\omega}$ as in Lemma \ref{LDMR}.  Define $\mathcal{G} : [Y \cap Z]^{\leq\omega} \to [\mathcal{CL}(Y)]^{\leq\omega}$ as $\mathcal{G}(A) = \{\{x\} : x \in A\} \cup \{\textnormal{cl}_X(U_n) \cap Y : n \in \N\} \in [\mathcal{CL}(Y)]^{\leq\omega}$ for each $A \in [Y \cap Z]^{\leq\omega}$. It is clear that $\mathcal{G}$ is $\omega$-monotone. For each $A \in [Y \cap Z]^{\leq\omega}$ let $s_A = r_{\mathcal{G}(A)} \restriction Z$ and $\mathcal{O}(A) = \{N \cap Z : N \in \mathcal{N}(\mathcal{G}(A))\}$. As above we can verify that the assignments $A \to \mathcal{O}(A)$ and $A \to r_{\mathcal{G}(A)}$ witness that $Y \cap Z$ is monotonically retractable in $Z$.
\endproof

The next notion, which is the relative version of monotone $\omega$-stability, will help us to characterize the Valdivia compact spaces.

\begin{definition} Let $Y$ be a subspace of $X$.
We say that $Y$ is \textit{monotonically $\omega$-stable} in $X$ if there exits an $\omega$-monotone map $\mathcal{N} : [C_p(X)]^{\leq\omega} \to [\mathcal{P}(Y)]^{\leq\omega}$ such that $\mathcal{N}(A)$ is a network of $\Delta_{\textnormal{cl}(A)}\restriction_Y$, for all $A \in [C_p(X)]^{\leq\omega}$.
\end{definition}

\begin{remark}\label{ENF}
Let $f: X \to Z$ be an onto map, let $\mathcal{N} \subseteq \mathcal{P}(X)$ and let $G \subseteq C_p(Z)$ be such that $\Delta_{G} : Z \to \R^G$ is an embedding. Since $\Delta_{f^\ast(G)} = \Delta_{G} \circ f$, then the family $\mathcal{N}$ is a network of $f$ if and only if it is a network of $\Delta_{f^\ast(G)}$.
\end{remark}

\begin{theorem}\label{PMRME}
Let $Y$ be a dense subspace of $X$. If $Y$ is monotonically retractable in $X$, then $Y$ is monotonically $\omega$-stable in $X$.
\end{theorem}

\proof
Assign to each $\mathcal{F} \in [\mathcal{CL}(Y)]^{\leq\omega}$ a retraction $r_\mathcal{F} : X \to Y$ and a family $\mathcal{N}(\mathcal{F}) \in [\mathcal{P}(Y)]^{\leq\omega}$ as in Lemma \ref{LDMR}. Define the map $\mathcal{G} : [C_p(X)]^{\leq\omega} \to [\mathcal{CL}(Y)]^{\leq\omega}$ by $\mathcal{G}(E) = \{f^{-1}(\textnormal{cl}(B)) \cap Y : f \in E \textnormal{ and } B \in \mathcal{B}(\R)\}$ for each $E \in [C_p(X)]^{\leq\omega}$. Note that $\mathcal{G}$ is $\omega$-monotone. Let $\mathcal{O} = \mathcal{N} \circ \mathcal{G} : [C_p(X)]^{\leq\omega} \to [\mathcal{P}(Y)]^{\leq\omega}$ which is evidently $\omega$-monotone. We will prove that for every $E \in [C_p(X)]^{\leq\omega}$ the family  $\mathcal{O}(E)$ is a network of $\Delta_{\textnormal{cl}(E)}\restriction_Y$. Fix $E \in [C_p(X)]^{\leq\omega}$.  Set $r = r_{\mathcal{G}(E)}$. Choose an arbitrary map $f \in E$. From property $(d)$ of Lemma \ref{LDMR}, we know that  $r \restriction_Y (F) \subseteq F$ for each $F \in \mathcal{G}(E)$. It follows from Lemma 4.15 in \cite{rjs-tka} that $f \restriction_Y = f \circ r \restriction_Y$. The density of $Y$ in $X$ guarantees that $f = f \circ r \in r^\ast(C_p(r(Y)))$. Since $f \in E$ was taken arbitrarily, we conclude that $E \subseteq r^\ast(C_p(r(Y)))$. Since every retraction map is $\R$-quotient, then $r^\ast(C_p(r(Y)))$ is closed in $C_p(X)$ and, hence, $\textnormal{cl}(E) \subseteq r^\ast(C_p(r(Y)))$. According to Lemma \ref{LDMR} $(b)$, we know that the family $\mathcal{O}(E)$ is a network of $r \restriction_Y$. Then the Remark \ref{ENF} implies that $\mathcal{O}(E)$ is a network of $\Delta_{(r \restriction_Y)^\ast(C_p(r(Y)))} = \Delta_{r^\ast(C_p(r(Y)))} \restriction_Y$. Therefore, $\mathcal{O}(E)$ is a network of $\Delta_{\textnormal{cl}(E)}\restriction_Y$.
\endproof

\section{Another characterization of Valdivia compact spaces}

Now we are ready to provide a characterization of Valdivia compact spaces using relative monotonically retractable spaces.

\begin{theorem}\label{DMRCRS}
Let $Y$ be a dense subspace of $X$. If $Y$ is (commutatively) monotonically retractable in $X$, then $Y$ is induced by an (a commutative) $r$-skeleton in $X$.
\end{theorem}

\proof
Assign to each $\mathcal{F} \in [\mathcal{CL}(Y)]^{\leq\omega}$ a retraction $r_\mathcal{F} : X \to Y$ a set $\mathcal{D}(\mathcal{F}) \in [Y]^{\leq\omega}$ and a family $\mathcal{N}(\mathcal{F}) \in [\mathcal{P}(Y)]^{\leq\omega}$ closed under finite intersections satisfying $(a)$ - $(d)$ of Lemma \ref{LDMR}. Set $\Gamma = [\mathcal{CL}(Y)]^{\leq \omega}$. We will prove that $\{r_\mathcal{F} : \mathcal{F} \in \Gamma\}$ is an $r$-skeleton in $X$. The conditions $(i)$ - $(iii)$ can be verified almost as in Theorem 4.3 in \cite{cas-rjs}.

To prove $(iv)$, fix $x \in X$. Let $U$ be an open neighborhood of $x$ in $X$. By regularity, we can find an open neighborhood $V$ of $x$ in $X$ such that $x \in V \subseteq \textnormal{cl}_X(V) \subseteq U$. Set $F = \textnormal{cl}_X(V) \cap Y$ and $\mathcal{F} = \{F\} \in \Gamma$. Choose $\mathcal{G} \in \Gamma$ such that $\mathcal{F} \subseteq \mathcal{G}$. We then have that $r_\mathcal{G}(F) \subseteq F$. The density of $Y$ in $X$ implies that $x \in \textnormal{cl}_X(F)$ and, hence, $r_\mathcal{G}(x) \in r_\mathcal{G}(\textnormal{cl}_X(F)) \subseteq \textnormal{cl}_Y (r_\mathcal{G}(F)) \subseteq F \subseteq \textnormal{cl}_X (V) \subseteq U$. Therefore, $x = \lim_{s \in \Gamma}r_\mathcal{F}(x)$.

Finally, it is easy to verify that $Y$ is induced by the $r$-skeleton $\{r_\mathcal{F} : \mathcal{F} \in \Gamma\}$ in $X$. In the commutative case, by Lemma \ref{LDMR},  the $r$-skeleton $\{r_\mathcal{F} : \mathcal{F} \in \Gamma\}$ is also commutative.
\endproof

Now we give conditions under which Theorem \ref{DMRCRS} has a converse.

\begin{theorem}\label{MWERSMR}
Let $Y$ be a subspace of $X$. If $Y$ is monotonically $\omega$-stable in $X$ and is induced by an (a commutative) $r$-skeleton in $X$, then $Y$ is (commutatively) monotonically retractable in $X$.
\end{theorem}

\proof Let  $\mathcal{O} : [C_p(X)]^{\leq\omega} \to [\mathcal{P}(Y)]^{\leq\omega}$ be an $\omega$-monotone map such that $\mathcal{O}(M)$ is a network of $\Delta_{\textnormal{cl}(M)}\restriction_Y$, for all $M \in [C_p(X)]^{\leq\omega}$. Let $\{r_s : s \in \Gamma\}$ be an $r$-skeleton on $X$ such that $Y = \bigcup\{r_s : s \in \Gamma\}$. For each $y \in Y$ fix $s_y \in \Gamma$ such that $r_{s_y}(y) = y$. Let $\gamma : [Y]^{\leq\omega} \to \Gamma$ be the map induced by the assignment $y \to s_y$ as in Lemma \ref{LAMG}.
Choose $F \in [Y]^{<\omega}$. The space $r_{\gamma(F)}(X)$ is cosmic, so the spaces $C_p(r_{\gamma(F)}(X))$ and $r_{\gamma(F)}^\ast(C_p(r_{\gamma(F)}(X)))$ are cosmic too. Fix a countable dense subset $D_F$ of $r_{\gamma(F)}^\ast(C_p(r_{\gamma(F)}(X)))$. Define a map $\mathcal{M} : [Y]^{\leq\omega} \to [C_p(X)]^{\leq\omega}$ by $\mathcal{M}(A) = \bigcup\{D_F : F \in [A]^{<\omega}\}$ for each $A \in [Y]^{\leq\omega}$. Note that $\mathcal{M}$ is $\omega$-monotone. For each $A \in [Y]^{\leq\omega}$, set $r_A := r_{\gamma(A)}$ and $\mathcal{N}(A) := \mathcal{O}(\mathcal{M}(A))$. We shall verify that the assignments $A \to r_A$ and $A \to \mathcal{N}(A)$ witness that $Y$ is monotonically retractable in $X$. It is clear that $\mathcal{N}$ is $\omega$-monotone. Given $A \in [Y]^{\leq\omega}$, to see that $A \subseteq r_A(X)$, choose $y \in A$. Then $s_y \leq \gamma(\{y\}) \leq \gamma(A)$ and $y = r_{s_y}(y) = r_{\gamma(A)}(r_{s_y}(y)) = r_{\gamma(A)}(y) = r_A(y) \in r_A(X)$. Finally we have the following.\medskip

\textbf{Claim.} If $A \in [Y]^{\leq\omega}$, then the family $\mathcal{N}(A)$ is a network of $r_{A} \restriction_Y$.\medskip

{\it Proof of the Claim.} Fix $A \in [Y]^{\leq\omega}$ and set $G := C_p(r_A(Y))$. By Remark \ref{ENF}, it is enough to show that the family $\mathcal{N}(A)$ is a network of $\Delta_{(r_{A}\restriction_Y)^\ast(G)}$. We know that the family $\mathcal{N}(A) = \mathcal{O}(\mathcal{M}(A))$ is a network of $\Delta_{\textnormal{cl}(\mathcal{M}(A))} \restriction_Y$. So we only need to verify that $(r_{A}\restriction_Y)^\ast(G) \subseteq \textnormal{cl}(\mathcal{M}(A)) \restriction_Y$. Fix $f \in (r_{A}\restriction_Y)^\ast(G)$. Then we can select $g \in G$ so that $f = g \circ (r_A \restriction_Y)$. To prove that $f \in \textnormal{cl}(\mathcal{M}(A)) \restriction_Y$  we consider an arbitrary canonical open neighborhood $W = [x_1,\ldots,x_k;U_1,\ldots,U_k]$ of $g \circ r_A$ in $C_p(X)$.  Let $V = [x_1,\ldots,x_k ; g^{-1}(U_1),\ldots,g^{-1}(U_k)]$ which is open in  $C_p(X,r_A(Y))$ and note that $r_A = r_{\gamma(A)} \in V$. Put $A=\bigcup\{F_n : n \in \N\}$ where $F_n \in [A]^{<\omega}$ and $F_n \subseteq F_{n+1}$ for each $n \in \N$. We know that $r_{\gamma(F_n)}(X) \subseteq r_{\gamma(A)}(X)$, for each $n \in \N$, and $r_{\gamma(A)}(x) = \lim_{n \to \infty} r_{\gamma(F_n)}(x)$, for every $x \in X$. Then we can choose $m \in \N$ such that $r_{\gamma(F_{m})} \in V$. It follows that $g \circ r_{\gamma(F_{m})} \in W$. Since $g \circ r_{\gamma(F_{m})} \in r_{\gamma(F_{m})}^\ast(C_p(r_{\gamma(F_{m})}(X)))$ and $D_{F_m}$ dense in $r_{\gamma(F_m)}^\ast(C_p(r_{\gamma(F_m)}(X)))$, we can find $h \in D_{F_m} \cap W \subseteq \mathcal{M}(A) \cap W$ which implies that $W \cap \mathcal{M}(A) \not= \emptyset$. Since $W$ was taken arbitrarily, we conclude that $g \circ r_A \in \textnormal{cl}(\mathcal{M}(A))$. Therefore, $f = (g \circ r_A) \restriction_Y \in \textnormal{cl}(\mathcal{M}(A))\restriction_Y$.\medskip

Assume that, in addition, the $r$-skeleton $\{r_s : s \in \Gamma\}$ which induces $Y$ is commutative. Then $r_A \circ r_B = r_B \circ r_A$ for all $A,B \in [Y]^{\leq\omega}$. Therefore,  $Y$ is commutatively monotonically retractable in $X$.
\endproof

The next result is a direct consequence of Theorems \ref{PMRME}, \ref{DMRCRS} and \ref{MWERSMR}.

\begin{corollary}\label{MRif and only ifMWERS}
Let $Y$ be a dense subspace of $X$. Then $Y$ is (commutatively) monotonically retractable in $X$ if and only if $Y$ is monotonically $\omega$-stable in $X$ and is induced by a (commutative) $r$-skeleton in $X$.
\end{corollary}

We know that a Compact space is Valdivia if and only if it admits a dense subset induced by a commutative $r$-skeleton (\cite{kbs-mich}). Moreover any compact space $X$ is monotonically $\omega$-stable  (\cite[Cor. 4.14]{rjs2}) and, as a consequence, any of its subspaces is monotonically $\omega$-stable in $X$. Hence, we obtain the next corollary.

\begin{corollary}\label{CVC}
A compact space $X$ is Valdivia if and only if it has a dense subset $Y$ which is monotonically retractable in $X$.
\end{corollary}

\section{$c$-skeletons}

Let us recall the the notion of $q$-skeleton which will be dualized in terms of monotone families of closed sets and open sets.

\begin{definition}\label{DQS}\cite{sal-rey} Let $X$ be a space. Consider a family $\{(q_s,D_s): s \in \Gamma\}$, where $q_s : X \to X_s$ is an $\R$-quotient map and $D_s$ is a countable subset of $X$ for each $s \in \Gamma$. We say that $\{(q_s,D_s): s \in \Gamma\}$ is a \textit{$q$-skeleton} on  $X$ if:
\begin{enumerate}[(i)]
	\item the set $q_s(D_s)$ is dense in $X_s$,
  \item if $s,t \in \Gamma$ and $s \leq t$, then there exists a continuous onto  map $p_{t,s} : X_t \to X_s$ such that $q_s = p_{t,s}\circ q_t$, and
  \item the assignment $s \to D_s$ is $\omega$-monotone.
\end{enumerate}
In addition, if $C_p(X) = \bigcup_{s \in \Gamma}q_s^\ast(C_p(X_s))$, then we say that the $q$-skeleton  is  {\it full}.
\end{definition}

In order to get a dual concept to $q$-skeleton, we introduce the following notion.

\begin{definition}\label{DCS} Let $X$ be a space and let $\{(F_s,\mathcal{B}_s): s \in \Gamma\}$ be a family of subsets of $X$ such that  $F_s$ is a closed subsets of $X$ and $\mathcal{B}_s$ is a countable family of open subset of $X$ for each $s \in \Gamma$. We say that $\{(F_s,\mathcal{B}_s): s \in \Gamma\}$ is a \textit{$c$-skeleton} on  $X$ if:
\begin{enumerate}[(i)]
	\item for each $s \in \Gamma$, $\mathcal{B}_s$ is a base for a topology denoted by  $\tau_s$ on $X$ and there exist a Tychonoff space $Z_s$ and a continuous map $g_s : (X,\tau_s) \to Z_s$ which separates the points of $F_s$,
  \item if $s,t \in \Gamma$ and $s \leq t$, then $F_s \subseteq F_t$, and
  \item the assignment $s \to \mathcal{B}_s$ is $\omega$-monotone.
\end{enumerate}
In addition, if $X = \bigcup_{s \in \Gamma}F_s$, then we say that the $c$-skeleton  is  {\it full}.
\end{definition}

Let us observe that the existence of a (full) $c$-skeleton in a space $X$ implies the existence of a (full) $c$-skeleton in any subspace of $X$.
This makes a big difference with the $r$-skeletons and $q$-skeletons.

\medskip

Now we establish the duality among the above concepts. First, we introduce some notation.
\medskip

Given $n \in \N \setminus \{0\}$, $N_1,\ldots,N_n \in \mathcal{P}(X)$ and $U_1,\ldots,U_n \in \mathcal{P}(\R)$, put
 $$[N_1,\ldots, N_n;U_1,\ldots, U_n] = \{f \in C_p(X) : \forall i \leq n \big(f(N_i) \subseteq U_i\big) \}.$$
 Consider the map $\mathcal{W} : [\mathcal{P}(X)]^{\leq\omega} \to [\mathcal{P}(C_p(X))]^{\leq\omega}$ defined by
 $$\mathcal{W}(\mathcal{N}) = \{[N_1,\ldots, N_n; B_1,\ldots, B_n] :  n \in \mathbb{N} \setminus \{0\} \ \text{and} \ \forall i \leq n  \big(N_i \in \mathcal{N} \wedge B_i \in \mathcal{B}(\R) \big) \},$$
 for each $\mathcal{N} \in [\mathcal{P}(X)]^{\leq\omega}$. It is straightforward to verify that $\mathcal{W}$ is $\omega$-monotone. For each $D \in [X]^{\leq\omega}$, let $\mathcal{W}_0(D) = \mathcal{W}(\{\{x\} : x \in D\})$ which is the family of all canonical open sets in $C_p(X)$ with support in $D$. It is easy to see that the map $\mathcal{W}_0 : [X]^{\leq\omega} \to [\mathcal{P}(C_p(X))]^{\leq\omega}$ is also $\omega$-monotone.

\begin{proposition}\label{DCQ}
If $X$ has a (full) $c$-skeleton, then $C_p(X)$ has a (full) $q$-skeleton.
\end{proposition}

\proof
Let $\{(F_s,\mathcal{B}_s): s \in \Gamma\}$ be a $c$-skeleton in $X$ and set $\mathcal{B} = \bigcup\{\mathcal{B}_s : s \in \Gamma\}$. For each nonempty set $N \in \mathcal{W}(\mathcal{B})$ fix a map $d_N \in N$. For each $s \in \Gamma$ we define $X_s := \pi_{F_s}(C_p(X))$, $q_s := \pi_{F_s} : C_p(X) \to X_s$ and  $D_s := \{d_N : N \in \mathcal{W}(\mathcal{B}_s)\}$. We shall verify that $\{(q_s,D_s): s \in \Gamma\}$ is a $q$-skeleton on $C_p(X)$.\medskip

(i) Fix $s \in \Gamma$. We shall verify that $q_s(D_s)$ is dense in $X_s$. Consider a nonempty canonical open set $U = [x_1,\ldots,x_n;I_1,\ldots,I_n]$ in $X_s$, where $x_1,\ldots,x_n \in F_s$ are pairwise distinct and $I_1,\ldots,I_n \in \mathcal{B}(\R)$. Since $g_s$ separates the points of $F_s$ and $Z_s$ is Tychonoff, we can find a continuous map $h \in C_p(Z_s)$ such that $h(g(x_i)) \in I_i$ for $i = 1,\ldots,n$. Let $f = h \circ g_s$ and note that $f(x_i) \in I_i$ for $i = 1,\ldots,n$. Since $f : (X,\tau_s) \to \R$ is continuous and $\mathcal{B}_s$ is a base for $\tau_s$ we can find $B_i \in \mathcal{B}_s$ such that $x_i \in B_i$ and $f(B_i) \subseteq I_i$ for $i = 1,\ldots,n$. If $N = [B_1,\ldots,B_n;I_1,\ldots,I_n]$, then $f \in N$ and $N \in \mathcal{W}(\mathcal{B}_s)$. Then $d_N \in D_s$. Note that $q_s(d_N)(x_i) = d_N \restriction_{ F_s}(x_i) = d_N(x_i) \in d_N(B_i) \subseteq I_i$ for $i = 1,\ldots,n$. Hence, $q_s(d_N) \in U$. Therefore, $q_s(D_s)$ is dense in $X_s$.\medskip

(ii) If $s,t \in \Gamma$ and $s \leq t$, then the map $p_{t,s} : X_t \to X_s$ given by $p_{t,s}(h) = h \restriction_{ F_s}$, for each $h \in X_t$, is continuous, onto and satisfies that $q_s = p_{t,s}\circ q_t$.\medskip

(iii) Since $\mathcal{W}_0$ is $\omega$-monotone, the assignment $s \to D_s$ also is $\omega$-monotone. \medskip

Assume now that  $X = \bigcup_{s \in \Gamma}F_s$. Let $\phi \in C_p(C_p(X))$. By applying Corollary 1.7.8. in \cite{ark-tko}, we can find a countable set $A  =\{a_n : n \in \N\} \subseteq X$ and a continuous map $\psi : \pi_A(C_p(X)) \to \R$ such that $\phi = \psi \circ \pi_A$. For each $n \in \N$ we can choose $s_n \in \Gamma$ such that $a_n \in F_n$ and $s_n \leq s_{n+1}$. Let $s = \sup\{s_n : n \in \N\}$. Notice that $A \subseteq F_s$. Consider the continuous map $\pi_{F_s,A} : X_s \to \pi_A(C_p(X))$ given by $\pi_{F_s,A}(h) = h \restriction_A$ for each $h \in X_s$. Then
$$
\phi = \psi \circ \pi_A = \psi \circ \pi_{F_s,A} \circ \pi_{F_s} = \psi \circ \pi_{F_s,A} \circ q_s = q_s^\ast(\psi \circ \pi_{F_s,A}) \in q_s^\ast(C_p(X_s)).
$$
Hence $C_p(C_p(X)) = \bigcup_{s \in \Gamma}q_s^\ast(C_p(X_s))$ and the $q$-skeleton $\{(q_s,D_s): s \in \Gamma\}$ in $C_p(X)$ is full.
\endproof

\begin{proposition}\label{DQC}
If $X$ has a (full) $q$-skeleton, then $C_p(X)$ has a (full) $c$-skeleton.
\end{proposition}

\proof
Let $\{(q_s,D_s): s \in \Gamma\}$ be a $q$-skeleton in $X$. Fix  $s \in \Gamma$. Since $q_s$ is a quotient map, the set $F_s := q_s^\ast(C_p(X_s))$ is closed in $C_p(X)$. On the other hand note that $\mathcal{B}_s = \mathcal{W}_0(D_s)$ is a countable family of canonical open subsets of $C_p(X)$. We shall prove that $\{(F_s,\mathcal{B}_s): s \in \Gamma\}$ is a $c$-skeleton on $C_p(X)$.\medskip

(i) Fix $s \in \Gamma$. It is clear that $\mathcal{B}_s$ is a base for a topology $\tau_s$ on $C_p(X)$. Consider the Tychonoff space $Z_s = \pi_{D_s}(C_p(X))$. Then $g_s = \pi_{D_s} : (C_p(X),\tau_s) \to Z_s$ is a continuous and onto map. We will verify that $g_s$ separates the points of $F_s$. Assume that $f_1,f_2 \in C_p(X)$ and $g_s(f_1) = g_s(f_2)$. Choose $h_1,h_2 \in C_p(X_s)$ such that $f_1 = h_1 \circ q_s$ and $f_2 = h_2 \circ q_s$. For each $x \in D_s$ we have that $h_1(q_s(x)) = f_1(x) = f_2(x) = h_2(q_s(x))$; that is, $h_1 \restriction_{q_s(D_s)} = h_2 \restriction_{q_s(D_s)}$. Since $q_s(D_s)$ is dense in $X_s$ we must have that $h_1 = h_2$, i.e. $f_1 = f_2$. This shows that $g_s$ separates the points of $F_s$.\medskip

(ii) Fix $s,t \in \Gamma$ with $s \leq t$. Then there exists a continuous onto  map $p_{t,s} : X_t \to X_s$ such that $q_s = p_{t,s}\circ q_t$. It follows that
$$
F_s = q_s^\ast(C_p(X_s)) = (p_{t,s}\circ q_t)^\ast(C_p(X_s)) = q_t^\ast \circ  p_{t,s}^\ast(C_p(X_s)) \subseteq q_t^\ast(C_p(X_t)) = F_t.
$$

(iii) It is clear that the assignment $s \to \mathcal{B}_s$ is $\omega$-monotone.\medskip

If the $q$-skeleton $\{(q_s,D_s): s \in \Gamma\}$ is full, then $C_p(X) = \bigcup_{s \in \Gamma}q_s^\ast(C_p(X_s)) = \bigcup_{s \in \Gamma}F_s$, and hence the $c$-skeleton $\{(F_s,\mathcal{B}_s): s \in \Gamma\}$ on $C_p(X)$ is full as well.
\endproof

It was proved in \cite{sal-rey} that: If $X$ has a full $q$-skeleton, then any countably compact subspace of $C_p(X)$ has a full $r$-skeleton. Besides, a compact space is Corson if and only if $C_p(X)$ has a full $q$-skeleton. By using these facts and the proposition from above, we have the following corollaries.

\begin{corollary}
If $X$ has a full $c$-skeleton and $Y$ is a countably compact subspace of some iterated function space $C_{p,2n}(X)$ for some $n \in \N$, then $Y$ has a full $r$-skeleton. In particular, if $Y$ is compact, then $Y$ is Corson.
\end{corollary}

\begin{corollary}
A compact space $X$ is Corson if and only if has a full $c$-skeleton.
\end{corollary}

Almost repeating the proof of Proposition 4.7 in \cite{sal-rey} we may prove the following proposition.

\begin{proposition}
If $X$ has a dense subspace $Y$ which is monotonically $\omega$-stable in $X$, then $X$ has a full $q$-skeleton.
\end{proposition}

We finish the paper with a list of unsolved questions.

\begin{question}
Let $X$ be a countably compact space, is it true $X$ has a full $c$-skeleton if and only if $X$ has a full $r$-skeleton?
\end{question}

By Propositions \ref{DCQ} and \ref{DQC} we know that if $X$ has a full $q$-skeleton, then $C_p(C_p(X))$ has a full $q$-skeleton. However, we do not know the answer to the following question.

\begin{question}
Must $X$ have a full $q$-skeleton whenever $C_p(C_p(X))$ has a full $q$-skeleton?
\end{question}

\end{document}